\documentclass[leqno]{amsart}

\usepackage{amscd,amssymb,amsmath}
\usepackage{mathabx}
\usepackage{float}
\usepackage{multirow}

\usepackage[utf8]{inputenc}
\usepackage{enumitem}

\usepackage{url}
\usepackage[OT1]{fontenc}
\usepackage[english,polish]{babel}
\usepackage{array,multirow,makecell}
\setcellgapes{2pt}
\makegapedcells
\newcolumntype{R}[1]{>{\raggedleft\arraybackslash }b{#1}}
\newcolumntype{L}[1]{>{\raggedright\arraybackslash }b{#1}}
\newcolumntype{C}[1]{>{\centering\arraybackslash }b{#1}}

\usepackage{amsmath,amsthm,amsfonts,latexsym,amscd,amssymb,enumerate}

\usepackage{amscd,amssymb,amsmath}

\usepackage{color}
\usepackage{colortbl}
\usepackage{tikz}
\usepackage{xcolor}
\usepackage{graphicx}
\usepackage{vntex}

\usepackage{amsmath,amsfonts,amsthm,amssymb}

\newtheorem{theorem}{Theorem}[section]
\newtheorem{lemma}[theorem]{Lemma}

\newtheorem{corollary}[theorem]{Corollary}
\newtheorem{proposition}[theorem]{{ Proposition}}

\newtheorem{definition}[theorem]{{Definition}}

\newtheorem{remark}[theorem]{{Remark}}

\newcommand\R{{\mathbb R}}

\newcommand\N{{\mathbb N}}
\newcommand\C{{\mathbb C}}
\newcommand\Z{{\mathbb Z}}

\newcommand{\U}{\mathcal{U}}

\newcommand{\D}{\mathcal{D}}

\def\strate{S}

 \def\codim{{\rm codim}}

\def\arho{d}
 
\def\jalpha{j}

\def\sangle{\sphericalangle}

\begin{document}


\title[Chern-Schwartz-MacPherson classes in the Lipschitz framework]
{Chern-Schwartz-MacPherson classes\\ in the point of view of Obstruction Theory \\ 
and Lipschitz framework.}


\medskip

\thanks{The first named author was supported by CNRS (UMR 7373) and Aix-Marseille University (I2M). 
The third named author was supported by the Brazilian S\~ao Paulo Research Foundation (FAPESP) 
grant \#2023/07802-7.}

\author{Jean-Paul Brasselet }
\address{I2M CNRS,  Aix-Marseille University, Marseille, France.
\newline jean-paul.brasselet@univ-amu.fr}

\author{Tadeusz  Mostowski}
\address{Uniwersytet Warszawski, Polska, 
\newline tmostows@mimuw.edu.pl}

\author{Th\h{u}y Nguy\~\ecircumflex n Th\d{i} B\'ich} 
\address{S\~ao Paulo State University, UNESP/Ibilce, S\~ao Jos\'e do Rio do Preto, Brazil, \quad
\newline
bich.thuy@unesp.br}
\maketitle \thispagestyle{empty}

\selectlanguage{english}

\begin{abstract} 
Since Chern and Grothendieck, Chern's characteristic class theory has made significant progress. 
In particular with regard to the classes of singular varieties. 
Conjectured by Grothendieck and Deligne and demonstrated by MacPherson, 
Chern classes of singular varieties have been defined in several ways, 
such as using polar varieties, Lagrangian theory...
However, the initial definition using obstruction theory, 
due to Marie-H\'el\`ene Schwartz, has been forgotten. 
Despite the simple ideas that enabled the obstruction definition, their implementation using 
Whitney stratifications requires delicate and technical  constructions. 
In the present article, we show that in the Lipschitz framework, 
the ideas of Marie-H\'el\`ene  Schwartz lead to a simplified definition and construction of Chern classes 
of complex analytic varieties.
\end{abstract}

\bigskip

\section{Introduction.} 

In his creative article \cite{Ch2}, Chern gives various definitions of  classes of hermitian manifolds.
One of them uses obstruction theory for the  construction of frames tangent to the manifold. 
In 1965,  M.-H. Schwartz provides a definition of 
Chern classes for singular analytic varieties, using obstruction theory in the framework of Whitney stratifications. 
Her construction  involves rather technical and delicate  systems 
of tubular neighbourhoods of the strata. 

The aim of this article is to show that the M.-H. Schwartz' ideas for defining Chern classes 
in the singular setting can be implemented in the framework of Lipschitz stratifications, 
instead of Whitney ones and that, in this framework, the M.-H. Schwartz' construction is simplified. 
\medskip

In this article, manifolds are smooth analytic complex varieties. 
A complex $r$-field tangent to a subset $A$ of a complex manifold $M$ is a set of $r$ continous 
complex vectors 
fields tangent to the manifold. When the vectors are $\C$-linearly independent at all points of $A$,
we will speak of an $r$-frame. A singular point of an $r$-field is a 
point at which the $r$ vectors fail to be $\C$-linearly independent.
The classical obstruction theory (see \cite{Steenrod1951}) says  that, on an $m$-dimensional complex manifold $M$, 
it is possible to construct such a complex $r$-field 
without singularities on the  $2m-2r +1$-skeleton  
of a  triangulation of $M$ and with isolated singular points situated at the barycenters 
of the $2p=2(m-r+1)$ simplices of the triangulation.
Let $I(v^{(r)}, a)$ be the index of the $r$-field 
$v^{(r)}=(v_1, \ldots , v_r)$ at the barycenter $a$ 
of the $2p$-simplex $\sigma$. 
The Chern class (in $H^{2p}(M,\Z))$ is the class of the cocycle whose value on the 
simplex  $\sigma$ is the index $I(v^{(r)}, a)$. 
\medskip

In this article, we  consider a singular 
complex analytic variety $X$ embedded in the complex 
 Euclidean space $\C^m$. Note that our result applies to situations which admit 
 Lipschitz stratifications 
(T. Mostowski \cite{Mostowski1985}).
 
  When using a triangulation compatible with a stratification of $X$, 
 the obstruction dimension for the construction of a complex $r$-frame tangent to a stratum
 depends on the  considered stratum. It is then clearly impossible to define a cocycle 
 in the same way than in the smooth case because the dimension of the cocycle 
 depends on the dimension of each stratum. 

Given a stratification of $X$, a stratification of $\C^m$ is obtained by adding 
the stratum $\C^m \setminus X$. Let $(K)$ be a 
locally finite triangulation of $\C^m$ compatible with the stratification of $X$, 
{\it i.e.} every open simplex lies in one and only one stratum. 
The first idea of Marie-H\'el\`ene Schwartz is to consider in $\C^m$ a barycentric  
cellular decomposition $(D)$ dual of the triangulation $(K)$. 
A dual cell $d$ which is an element of $(D)$ is then topologically transverse to all strata $S$, 
that means (see Lemma \ref{delemme}) that 
$$ {\text {\rm codim}} (d) +  {\text {\rm codim}}  (S) =  {\text {\rm codim}} (d \cap S).$$
The intersection of a $2(m-r+1)$-dimensional  dual cell  with  
a stratum is a cell whose dimension is equal to the obstruction dimension for 
constructing an $r$-frame tangent to that stratum. It is then possible to define 
a $2(m-r+1)$-cocycle by attaching to each $2(m-r+1)$ cell of $(D)$ meeting $X$ 
the index of a stratified complex $r$-field  with isolated singularities
at the barycenters of the cells. Here, a stratified $r$-field 
means an $r$-field such that at each point $x$ of the stratification of $\C^m$, the $r$-field is 
tangent to the stratum containing $x$. 

The second idea of M.-H. Schwartz concerns the frames. 
When proving the Poincar\'e-Hopf theorem for singular varieties, M.-H. Schwartz has showed
that  considering stratified vector fields is not sufficient. She provides counter-examples
for the Poincar\'e-Hopf theorem using vector fields which are only stratified. 
She has showed that one has to consider radial extension of stratified vector fields.
The same situation appears  in the definition of characteristic classes (Chern classes). 
M.-H. Schwartz considers stratified $r$-fields,
obtained by suitable radial extensions. 
Such an $r$-field defines an 
element of the relative cohomology $H^{2(m-r+1)}(\C^m, \C^m\setminus X)$, called the 
Schwartz (cohomology) class of $X$. It is independent on the choices of stratification, triangulation and
of the radial extension. 

The M.-H. Schwartz' ideas are very clear and adapted. 
However, their implementation is rather technical. 
Under Whitney stratifications it uses  delicate constructions 
of systems of adapted tubular neighbourhoods. Moreover,  it needs to consider 
technical homotopies between frames issued from different strata. 
\medskip

Let us recall that, in his original article \cite{Mostowski1985}, T. Mostowski proved that analytic 
complex varieties admit Lipschitz stratifications. 

In this article, we use the M.-H. Schwartz' ideas in the framework of Lipschitz stratifications
instead of Whitney stratifications. That allows to simplify her construction. 
For example, the continuity of  $r$-frames comes directly from Lipschitz properties
without using homotopy between frames.

By Alexander duality isomorphism 
$$H^{2(m-r+1)}(\C^m, \C^m \setminus X) \buildrel{\cong}\over \longrightarrow H_{2(r-1)}(X),$$
the obtained homology classes are the MacPherson classes of the constructible 
function ${\bf 1}_X$  \cite{BS}. These homology classes have been defined by 
MacPherson  \cite{MacPherson1974}, answering 
a conjecture by Deligne and Grothendieck \cite[Note $87_1$]{Grothendieck2022}.

Nowadays, these classes are called (Chern-)Schwartz-MacPherson classes.
\medskip 

The paper is organized as follows: In the first section, 
we recall the classical obstruction theory
for manifolds and we explain the M.-H. Schwartz' idea for considering obstruction 
theory in the case of stratified 
singular varieties. The second section deals with Lipschitz stratifications: their definition and properties. 
In the third section, we show how Lipschitz properties allow to construct radial extensions 
of Lipschitz vector fields and frames. The last section is devoted to the definition 
of obstruction cocycles and classes in our setting.

\section{A few on obstruction theory}

\subsection{The smooth case}

Let $M$ be a locally compact complex analytic manifold with complex dimension $m$. 
The complex tangent bundle, denoted by $TM$ is locally trivial.

For a fixed  $r$, with $1\le r \le m$, an \emph{$r$-field} on a subset $A$ of $M$ 
is a set $v^{(r)}=\{v_1,\dots,v_r\}$
of $r$ continuous complex vector fields defined on $A$ and tangent to $M$. A singular point
of $v^{(r)}$ is a point where the vectors $\{v_1,\dots,v_r\}$ fail to be $\C$-linearly
independent. An \emph{$r$-frame} is a non-singular $r$-field.

The objective of the obstruction theory is to evaluate the obstruction to the construction of $r$
sections of $TM$ linearly independent  at each point, \emph{i.e.} an $r$-frame.

Let us denote by $V_r(TM)$ the fibre bundle over $M$
associated to $TM$ and whose fibre at the point $x$ of $M$ is the set of $r$-frames of 
the fiber $T_x M$ of $TM$ at $x$.
An $r$-frame  on a subset $A$ of $M$ is then a section of  $V_r(TM)$ over $A$. 
The typical fibre of $V_r(TM)$ is the Stiefel
manifold of complex $r$-frames in ${\C}^m$ denoted by $V_{r}(\mathbb{C}^m)$. 
In the same way that $TM$ is locally trivial, the associated bundle $V_r(TM)$ 
is locally trivial with the same triviality open subsets.

In this section, we use a triangulation of the manifold $M$. 
In the singular case (section \ref{tridec}), we will see the usefulness of 
considerning a cellular decomposition, with a similar construction.

We assume that the  triangulation $(K)$  of $M$ is sufficiently small so that every simplex $\sigma$ 
is included in an open subset  $U$ over which the bundle $V_r(TM)$ is trivial. 

Let us assume that one has an $r$-frame $v^{(r)}$ defined on the boundary 
$\partial \sigma$ of a $k$-simplex $\sigma$. It provides a map
\begin{equation}\label{demesb}
 \partial \sigma \buildrel {v^{(r)}}\over \longrightarrow V_r(TM)\vert_U 
\cong U\times V_{r}(\mathbb{C}^m) \buildrel {pr_2}\over \longrightarrow 
V_{r}(\mathbb{C}^m),
\end{equation}
where $pr_2$ is the second projection.
One obtains a map 
\begin{equation}\label{compoto}
\xi :  \mathbb{S}^{k-1} \cong \partial \sigma  \buildrel {pr_2 \circ v^{(r)}}\over 
\longrightarrow V_{r}(\mathbb{C}^m)
\end{equation} 
that defines an element of 
the homotopy group $\pi_{k-1} (V_{r}(\mathbb{C}^m))$, denoted by 
$[\xi]$. 

By classical homotopy theory, if  $[\xi]=0$, then the map
$  \partial \sigma \rightarrow V_{r}(\mathbb{C}^m)$ defined on the boundary 
of the simplex $\sigma$   can be extended inside $\sigma$.
In other words, if $[\xi]=0$, then 
 there is no obstruction for the extension of the 
 $r$-frame $v^{(r)}$, section of $V_r(TM)$, 
inside $\sigma$. This happens in particular when $\pi_{k-1} (V_{r}(\mathbb{C}^m)) = 0$.  

The homotopy groups $\pi_{*} (V_{r}(\mathbb{C}^m))$ are computed by
Stiefel and by Whitney (see \cite{Steenrod1951}). One has: 
\begin{equation}
\pi_{k-1} (V_{r}(\mathbb{C}^m)) =  
\begin{cases} 
0 &  \text{for $k< 2m-2r+2$}\cr
\mathbb{Z} &  \text{for $k= 2m-2r+2$}.\cr
\end{cases}
\label{sticomp}
\end{equation}
\goodbreak 

Moreover, one has the following results:

\begin{proposition}  \label{propcomp}
Let us consider  a complex $r$-frame $v^{(r)}$ defined on the boundary $\partial \sigma$ 
of the $k$-simplex $\sigma$. 
\begin{itemize}
\item  If $k<2(m-r+1)$,  then $[\xi]=0$ and one can extend 
the $r$-frame $v^{(r)}$ inside $\sigma$ without singularity.
\item  If  $k=2(m-r+1)$, then  one can extend 
the $r$-frame $v^{(r)}$  inside $\sigma$ with an isolated singularity at the barycenter 
$\widehat \sigma$ of the simplex. In that case, one obtains an index 
 $I(v^{(r)},\widehat \sigma)$ 
  defined by the class $[\xi]
 \in \pi_{k-1} (V_{r}(\mathbb{C}^m)) \cong \mathbb{Z}$. 
 The index measures the obstruction to the extension of 
$v^{(r)}$, defined on the boundary $\partial \sigma$, inside $\sigma$.
\end{itemize}
\end{proposition}
The dimension $2p= 2(m-r+1)$ is called the \emph{obstruction dimension} for the 
construction of a  complex $r$-frame tangent to $M$.

For $q=0,\ldots,2m$, we denote by $(K)^{(q)}$ the $(q)$-skeleton  
of $(K)$, union of simplices  
of dimension less than or equal to $q$. 

By Proposition \ref{propcomp}, 
it is possible to construct an $r$-frame $v^{(r)}$ without obstruction on the $(2p-1)$-skeleton 
$(K)^{(2p-1)}$
and with isolated singularities at the barycenters $\widehat \sigma$
of the $2p$-simplices $\sigma$ with indices $I(v^{(r)},\widehat \sigma)$.
The cochain $c^p : (K)^{(2p)} \to \Z$ defined by $c^p(\sigma) = I(v^{(r)},\widehat \sigma)$  
 is a $2p$-cocycle of $(K)$  (for instance see \cite{DK}). 
Its class in $H^{2p}(K ; \Z) = H^{2p}(M ; \Z) $,
denoted by $c^p(M)$,  is called the $p$-th cohomology Chern class of $M$.

The Chern class does not depend on the choices of the triangulation and of the $r$-frame. 

\subsection{The singular case} 

In the following, $X\subset \C^m$  will denote an analytic  
complex variety of complex dimension $n$, subset of a  compact region.


\subsubsection{Stratifications} 
 
A  {stratification} $({\mathcal S})$ of  $X$  is a family of 
 closed analytic subsets of $X$ 
$$
({\mathcal S}) : \qquad
X=X^n \supset X^{n-1} \supset X^{n-2} \supset \cdots  \supset X^1 
\supset X^0 \supset \emptyset = X^{-1},   
$$
where each stratum ${S}^j= X^j-X^{j-1}$ is either empty or a 
(non necessarily connected) smooth complex manifold of pure (complex) dimension $j$. 

A stratification of $X$ gives rise to a stratification of $\C^m$ by adding the stratum 
$\C^m \setminus X$. 

There are many ways to consider and define a stratification of an analytic variety. 
In her original work on characteristic classes for singular varieties,
Marie-H\'el\`ene Schwartz uses Whitney stratifications.
In this article, we consider Lipschitz stratifications ({\it cf.} Section \ref{Lip1}). 
For a complete presentation of stratifications, see for instance \cite{Trotman2022}.

\subsubsection{Triangulations and cellular decompositions}\label{tridec}

Let $(K)$ be a triangulation of $\C^m$ compatible with 
the stratification of $\C^m$. That means that every open simplex of $(K)$ lies in one
and only one stratum. The existence of such a triangulation has been proved by
Lefschetz \cite{Lefschetz1930} in the analytic case.

We consider a barycentric subdivision  $(\widehat K)$ of $(K)$. The barycenter 
of a simplex $(\sigma  \in K)$ will be denoted by $\widehat\sigma$. 
Every simplex in $\widehat K$ can be written as 
$$(\widehat \sigma_{i_1}, \widehat \sigma_{i_2}, \ldots , \widehat \sigma_{i_p})$$
where $\sigma_{i_1} < \sigma_{i_2} < \cdots < \sigma_{i_p}$. 
Here the symbol $\sigma < \sigma'$ means that the 
simplex $\sigma$ is a face of $\sigma'$. 

The dual cell of a simplex $\sigma$, denoted by $d(\sigma)$,  
 is the set of all (closed) simplices 
$\tau$ in $(\widehat K)$ such that $\tau \cap \sigma =
\{ \widehat\sigma \}$. That is the set of all simplices in $(\widehat K)$ on the form 
$(\widehat \sigma, \widehat \sigma_{i_1}, \ldots , \widehat \sigma_{i_k})$ 
with $\sigma < \sigma_{i_1} < \cdots < \sigma_{i_k}$.  

\begin{figure}[H]
\centering
\begin{tikzpicture}[scale=0.3]

\coordinate (A) at (0,0);
\coordinate (B) at (10,-3);
\coordinate (C) at (19.5,0.5);
\coordinate (D) at (2,7);
\coordinate (E) at (9,4);
\coordinate (F) at (16,8);
\coordinate (G) at (28,7);
\coordinate (H) at (8,14);
\coordinate (J) at (21,14);

\coordinate (bA) at (-4,0);
\coordinate (Gd) at (32,7);

\coordinate (K) at (6,0);
\coordinate (L) at (12,1);
\coordinate (M) at (4,4);
\coordinate (N) at (16,4);
\coordinate (P) at (21,5);
\coordinate (Q) at (6,8);
\coordinate (R) at (11,8);
\coordinate (S) at (16,12);
\coordinate (T) at (22,10);

\coordinate (AE) at (intersection of A--E and  5.5,0 -- 5.5,5);
\coordinate (AD) at (intersection of A--D and  0,4.2 -- 5,4.2);
\coordinate (AB) at (intersection of A--B and  5.5,0 -- 5.5,5);
\coordinate (BE) at (intersection of B--E and  0,1.3 -- 5,1.3);
\coordinate (BC) at (intersection of B--C and  15,0 -- 15,5);
\coordinate (CE) at (intersection of C--E and  14.5,0 -- 14.5,5);
\coordinate (FE) at (intersection of F--E and  12,0 -- 12,5);
\coordinate (HE) at (intersection of H--E and  0,7.5 -- 5,7.5);
\coordinate (DE) at (intersection of D--E and  5.5,0 -- 5.5,5);
\coordinate (DH) at (intersection of D--H and  4.5,0 -- 4.5,5);
\coordinate (FH) at (intersection of F--H and  12,0 -- 12,5);
\coordinate (FJ) at (intersection of F--J and  19,0 -- 19,5);
\coordinate (FG) at (intersection of F--G and  21,0 -- 21,5);
\coordinate (FC) at (intersection of F--C and  18,0 -- 18,5);
\coordinate (CG) at (intersection of C--G and  24,0 -- 24,5);
\coordinate (GJ) at (intersection of G--J and  24,0 -- 24,5);
\coordinate (HJ) at (intersection of H--J and  15,0 -- 15,5);

\coordinate (CZ) at (28,4);
\coordinate (GZ) at (27,11.5);
\coordinate (AZ) at (-2,4.5);
\coordinate (BZ) at (1,-2.5);

\fill[gray!45] (P)--(FG)--(T)--(FJ)--(S)--(FH)--(R)--(FE)--(N)--(FC)--(P);
\fill[gray!65] (L)--(CE)--(N)--(FE)--(R)--(HE)--(Q)--(DE)--(M)--(AE)--(K)--(BE)-- (L);

\fill[gray!45] (P) -- (FG)-- (T) -- (GJ) -- (G) -- (CG) -- (P) ;
\fill[gray!45] (A) -- (AD)--  (M) --(AE) --  (K) -- (AB) -- (A) ;

\draw[ultra thick,densely dashed] (bA) -- (A)--(E) -- (F) --(G) -- (Gd);

\draw[thick] (A)--(D)--(E);
\draw[thick] (A)--(B)--(E) -- (C) -- (B);
\draw[thick] (G)--(C)--(F);
\draw[thick] (E)--(H)--(D);
\draw[thick] (F)-- (H)--(J) -- (F);
\draw[thick] (G)--(J);

\draw[thick,dotted] (A)--(K)--(E)--(M)--(A);
\draw[thick,dotted] (E)--(L)--(C)--(N)--(E);
\draw[thick,dotted] (E)--(R)--(H)--(Q)--(E);
\draw[thick,dotted] (F)--(P)--(G)--(T)--(F);
\draw[thick,dotted] (N)--(F)--(R);
\draw[thick,dotted] (Q)--(D)--(M);
\draw[thick,dotted] (K)--(B)--(L);
\draw[thick,dotted] (P)--(C);
\draw[thick,dotted] (S)--(H);
\draw[thick,dotted] (T)--(J)--(S)--(F);

\draw (M)--(AE)--(K);
\draw (M)--(AD);
\draw (AB)--(K)--(BE)--(L)--(BC);
\draw (L)-- (CE)--(N)--(FE)--(R)--(HE)--(Q)--(DE)--(M);
\draw (DH)--(Q);
\draw (N)--(FC)--(P)--(FG)--(T)--(FJ)--(S)--(FH)--(R);
\draw (CG)--(P);
\draw (GJ)--(T);
\draw (HJ)--(S);

\coordinate (AZ) at (-2,4.5);
\coordinate (BZ) at (1,-2.5);

\draw [double](R)--(FE)--(N);

\draw[thick,dashed] (CG) -- (CZ);
\draw[thick,dashed] (GJ) -- (GZ);
\draw[thick,dashed] (AD) -- (AZ);
\draw[thick,dashed] (AB) -- (BZ);

\node at (N) {$\bullet$};
\node at (N)[below] {$B'$};
\node at (11,7.5) [above left]{$A'$};
\node at (R) {$\bullet$};
\node at (H) [above]{$C$};
\node at (FE) {$\bullet$};
\node at (11,10.15) {$\sigma_2$};
\node at (FE)[above] {$\widehat{\sigma_1}$};

\node at (-4,0) [above]{$\strate^j$};
\node at (9,3.6) [right]{A};
\node at (F)[below] {B};
\node at (14.5,7)[below] {$\sigma_1$};
\node at (7,6.8)[below] {$d(A)$};
\node at (13.4,5.4)[below] {$d(\sigma_1)$};
\node at (16,11)[below] {$d(B)$};

  \end{tikzpicture}
  \caption{Dual cells: The barycenter of $\sigma_0 = \{A \}$ is $A$ itself. The dual cell 
  of $\sigma_0$ is the  dark gray cell. The dual cell of $\sigma_1=AB$ is composed of the 
  two segments (double lines) $A' \widehat{\sigma_1}$ and $\widehat{\sigma_1} B'$. 
  The dual cell of the triangle $\sigma_2= ABC $ is the barycenter of $\sigma_2$ = $A'$.
The stratum $S^j$ is the thick dashed line} \label{dualcells2}
  \end{figure}
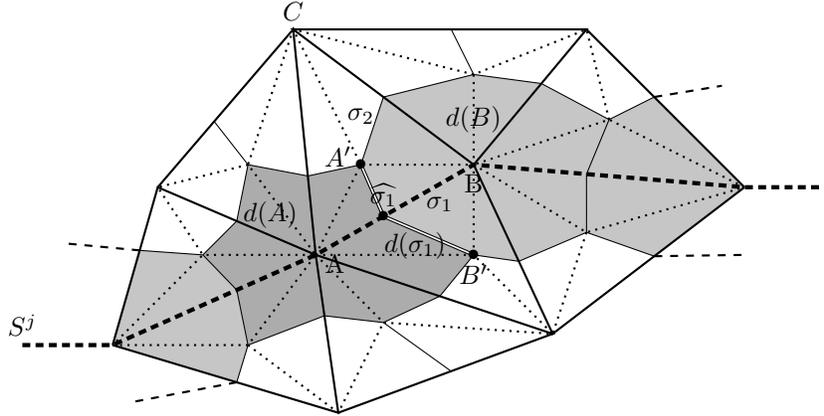

  The dual cells satisfy the nice properties (see for instance \cite{Munkres84}):
\begin{lemma}
{\rm (i)} The dual cell of a $k$-simplex is a $(2m-k)$-cell, homeomorphic to the unit ball 
${\mathbb B}^{2m-k} \subset \R^{2m-k}$ and its boundary is homeomorphic to 
the corresponding sphere ${\mathbb{S}}^{2m-k-1}$. \\
{\rm (ii)}
Given a triangulation $(K)$ of $\C^m$ compatible with 
the stratification of $\C^m$, 
the set of dual cells provides a cellular decomposition of $\C^m$, denoted by $(D)$ and called 
dual cellular decomposition associated to the barycentric subdivision $(\widehat K)$ of $(K)$. 
\end{lemma}

If the cell $d(\sigma)$ meets $X$, 
 the simplex $\sigma$ of which it is dual lies in $X$. 
The unique intersection point $\widehat\sigma = d(\sigma) \cap \sigma$ is the barycenter 
of $\sigma$ and is also the barycenter of $d(\sigma)$. It will be denoted by $\widehat d$.

The following straightforward Lemma plays an important role in the 
construction of M.-H. Schwartz \cite{Schwartz2000}.

\begin{lemma}\label{delemme} 
Let $({\mathcal S})$ be a stratification of $\C^m$ and   $(K)$ be a triangulation of $\C^m$ 
compatible with the stratification. Let $(D)$ be a cellular decomposition of $\C^m$ dual of $(K)$. 
Then the dual cells are (topologically) transverse to the strata.
\end{lemma}

Here, the topological transversality means that if ${S}^j$ is a 
stratum and $d$ is a dual cell, then the (real) codimensions in $\C^m$ satisfy 
$$\codim ( d \cap {S}^j) = \codim (d) + \codim ({S}^j).$$

\begin{remark}
The dual cells of all simplices $\sigma$ in the stratum ${S}^j$ form a cellular neighbourhood of 
${S}^j$ (see the grey neighbourhood of $S^j$ in Figure \ref{fig:dualcells2}). 
This fact will be useful in the following.
\end{remark} 

\subsubsection{Obstruction in a stratified space}\label{dobstruction}

According to Proposition \ref{propcomp}, 
the obstruction for constructing a stratified $r$-frame 
tangent to the stratum $\C^m \setminus X$ lies in dimension $2(m-r+1)$.
The obstruction for constructing 
an $r$-frame tangent to a stratum ${S}^j$ with complex dimension $j$ lies  in 
(real) dimension $2(j-r+1)$. 

 When using a triangulation $(K)$ of $\C^m$ compatible with a stratification of $\C^m$, 
 the obstruction dimension for the construction of a complex $r$-frame tangent to a stratum
 depends on  the  considered stratum.
 It is then clearly impossible to define a cocycle  in the same way than in the smooth case, 
 the dimension of the obstruction cocycle depends on the dimension of each stratum.

The first idea of Marie-H\'el\`ene Schwartz is to consider in $\C^m$ a  
cellular decomposition $(D)$ dual of $(K)$. In that case, as $(K)$ is compatible with the stratification,
the cells of $(D)$ are topologically transverse to the strata (Lemma \ref{delemme}). 
When not empty, the intersection of an $2(m-r+1)$ cell of $(D)$ with  a stratum 
${S}^j$ of (real) dimension $2 j$ 
is a cell of (real) dimension $2(j-r+1)$. 
That is exactly  the obstruction dimension for the 
construction of an $r$-frame tangent to ${S}^j$. It is then possible to define 
a $2p=2(m-r+1)$-cocycle $\gamma^p$ (see  formula (\ref{classeS})) 
by attaching to each $2(m-r+1)$ cell of $(D)$ meeting $X$ 
the index of a stratified complex $r$-field  at isolated singularities. 
We will use this first idea in the same way.

In the previous paper  (\cite[Section 2]{BMT}), we explicited the second idea of Marie-H\'el\`ene Schwartz.
To prove the Poincar\'e-Hopf theorem 
in the case of real semi-analytic varieties we need to use 
radial vector fields, which are obtained as a sum of parallel and transversal extensions. 
The same idea is valid for constructing the Schwartz classes.  
In particular this guarantees that the index of the considered  $r$-field 
 at an isolated singular point is the same, whether it is calculated on the stratum 
 or in the ambient space. Then, the cocycle $\gamma^p$ is well defined. 
In section \ref{ideas} 
we will  apply the M.-H. Schwartz'  idea in the context of Lipschitz stratifications.
The Schwartz class  $c^p(X)$ will be defined as the class of the cocycle $\gamma^p$.

\section{Lipschitz stratifications} \label{Lip1}

Motivated by a question of Sullivan (\cite{Siebenmann1979}, 1979) (for a complete 
discussion see \cite{Parusinski1994}), 
Lipschitz stratifications are defined  in 1985 by Tadeusz Mostowski \cite{Mostowski1985}. 
Mostowski proved existence of stratifications satisfying the so-called Lipschitz $L$-conditions for complex 
analytic varieties.  Some basic properties of Lipschitz stratifications 
are improved and clarified by Adam Parusi\'nski \cite{Parusinski1993}.

The first ingredient for the definition of Lipschitz stratifications is the notion of chains 
(see \cite{Mostowski1985,Mostowski1988,Parusinski1994}). 

\begin{definition}
[Definition of chains] \label{dechaine}
Let us fix a (real) constant  $K>1$. 
A {\sl chain} starting at a point  $q_{j_1} $  in an $j_1$-dimensional stratum  ${S}^{j_1}$ 
is a sequence of points $\{ q_{j_1} , \ldots, q_{j_t} , \ldots, q_{j_\ell} \}$
where $q_{j_k} \in {S}^{j_k},$ such that 
$$ X^{j_1} \supset X^{j_2}
\supset \cdots  \supset {X}^{j_t}  \supset \cdots  \supset {X}^{j_k}  \supset \cdots\supset {X}^{j_{\ell}} $$
(with complex dimensions $j_1 > j_2 > \cdots > j_t > \cdots >  j_\ell$) 
and where $j_t$ (depending on the constant $K$)  is the greatest integer for which 
\begin{equation}\label{dechaine1}
d(q_{j_1},{X}^{j_k}) \ge 2 K^2 \,  d(q_{j_1},{X}^{j_t}) \qquad
\text{ for all $k$ such that }j_t> j_k \ge j_ \ell
 \end{equation}
 and
\begin{equation}\label{dechaine2}  
\vert q_{j_1} - q_{j_t} \vert \le K \;  d(q_{j_1},{X}^{j_t}).
\end{equation}
\end{definition}
\bigskip

\begin{remark} {\rm 
The condition (\ref{dechaine1}) can be written 
$$d(q_{j_1},{X}^{j_t}) \le \frac{1}{2 K^2} d(q_{j_1} ,{X}^{j_k})$$
and means that when $q_{j_1}$ is going 
far away from ${X}^{j_t}$ to a stratum $S^{j_k}$ 
of lower dimension $j_k < j_t$, the distances 
$$d(q_{j_1},{X}^{j_t}) \text { and } d(q_{j_1},{X}^{j_k}) $$
are controlled by the inequality (\ref{dechaine1}).
}
\end{remark}

\begin{remark} {\rm 
The condition (\ref{dechaine2}) says that in the same stratum ${S}^{j_t}$, the distances 
 from the starting point $q_{j_1}$ to the point $q_{j_t}$ and 
 to the stratum ${S}^{j_t}$
are controlled by the inequality  (\ref{dechaine2}).
}
\end{remark}

The following example (here in $\R^3$, due to limitation of dimension in figures) 
illustrates the definition. 

\begin{figure}[H] \centering
\begin{tikzpicture}[scale=0.4] [domain=1:0.2]

\draw [thick,dashed] plot[domain=-1.65:0,samples=100] (2*\x,{(0.6*sin(\x r)+4)});
\draw [thick] plot[domain=0:1.7,samples=100] (2*\x,{(0.6*sin(\x r)+4)});
\draw [thick]plot[domain=-1.535:0,samples=100] (3*\x,{(-0.6*sin(\x r)+4)});
\draw [thick]plot[domain=1.15:1.5,samples=100] (3*\x,{(-0.6*sin(\x r)+4)});
\draw [thick,dashed] plot[domain=0:1.15,samples=100] (3*\x,{(-0.6*sin(\x r)+4)});

\draw [thick] plot[domain=-1.53:1.15,samples=100] (3*\x,{(0.6*sin(\x r)-4)});
\draw [thick]plot[domain=-1.15:1.52,samples=100] (3*\x,{(-0.6*sin(\x r)-4)});

\draw [thick]plot[domain=-1.15:1.15,samples=100] ({(0.6*sin(\x r)+4)},3*\x);
\draw [thick]plot[domain=-1.53:1.53,samples=100] ({(-0.6*sin(\x r)+4)},3*\x);

\draw  [thick,dashed]plot[domain=0:1.66,samples=100] ({(0.6*sin(\x r)-4)},2*\x);
\draw [thick]plot[domain=-1.53:0,samples=100] ({(0.6*sin(\x r)-4)},3*\x);

\draw [thick]plot[domain=0:1.525,samples=100] ({(-0.6*sin(\x r)-4)},3*\x);
\draw [thick]plot[domain=-1.16:0,samples=100] ({(-0.6*sin(\x r)-4)},3*\x);

    \fill [black,opacity=.5] (1.8,0.6) circle (4pt);
\node at (1,1.1) [right] {\scriptsize $q_{\jalpha_1}$};   
    \fill [black,opacity=.5] (1.8,0) circle (4pt);
\node at (1.8,0) [below] {\scriptsize $q_{\jalpha_2}$};   
    \fill [black,opacity=.5] (0,0) circle (4pt);
\node at (-0.5,0) [below] {\scriptsize $q_{\jalpha_3}$};   

\node at (2.8,4.6) [below] {\scriptsize $\strate^{\jalpha_1}$};   
\node at (3.3,0.2) [below] {\scriptsize $\strate^{\jalpha_2}$};   
\node at (-0.5,0) [above] {\scriptsize $\strate^{\jalpha_3}$};   


\draw (-4,0) -- (4,0);
\draw (0,-4) -- (0,4);
\node at (0,-5)  {.};
\end{tikzpicture}
\caption{Example of chain: $\jalpha_1= 2, \jalpha_t= \jalpha_2 = 1, \jalpha_3=0$, 
the chain is: $(q_{\jalpha_1}, \;q_{\jalpha_2},\;q_{\jalpha_3})$. 
}
\end{figure}
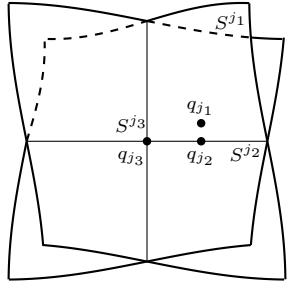

Before providing  the definition of Lipschitz stratifications, 
we fix some more notations: For  $q \in {S}^{j},$  let  
$$P_q: \C^m \to T_q {S}^{j} \subset T_q\C^m $$ 
be the orthogonal projection
onto the tangent space to the stratum $S^j$ at the point $q$ 
and  let $P_q^\perp = Id -P_q$ be the orthogonal projection onto the normal space 
$T_q^\perp  {S}^{j}$.

Identifying $T_q\C^m$ by $\C^m$, then $P_q$ and $P_q^\perp$ are linear maps from $\C^m$
to $\C^m$. The norm of the linear map $P_q$ can be defined as any norm of endomorphims of $\C^m$:
$$\Vert P_q\Vert := \sup_{\Vert x \Vert =1} \vert P_q(x) \vert.$$
The same stands for $P_q^\perp$. 

\begin{remark}{\rm 
An important remark (philosophy) 
concerning exponents and multiplicative constants 
is that exponents are important and interesting while 
multiplicative constants appearing in plenty of inequalities are not important.
We shall use occasionally the following notations:
$$f \lesssim g \Longleftrightarrow \text{ for some constant } C>0, \; f\le C g,$$
$$ f \sim g \Longleftrightarrow \vert f \vert  \lesssim  \vert g \vert  \lesssim \vert f \vert.$$}
\end{remark}

\begin{definition} [Lipschitz stratification] \label{deflip}
The stratification $(\mathcal S)$ is said a Lipschitz stratification if there exists a constant $C>0$ 
(called Lipschitz constant) such that \\
\noindent {a)} for every chain 
$ \{ q_{j_1},  q_{j_2}, \ldots,  q_{j_\ell} \}$ (Definition \ref{dechaine}), we have 
\begin{equation}\label{dequa1}
\Vert P^\perp_{q_{j_1}} P_{q_{j_2}} \cdots P_{q_{j_k}} \Vert \le C 
 \frac{\vert q_{j_1} - q_{j_2}\vert }{d(q_{j_1},X^{j_k -1})},
 \end{equation}
 for all  $1 < k \le \ell$. \\
\noindent {b)} If, further, 
$q_{j_1}$ and $q'_{j_1}$ are points in  ${S}^{j_1}$ such that 
$$\vert q'_{j_1}-q_{j_1} \vert \le \frac{1}{2K}\; d(q_{j_1},X^{j_1 -1})$$ 
(here $K$ is the constant in Definition \ref{dechaine}),  
then
\begin{equation}\label{dequa2}
\Vert (P_{q'_{j_1}} - P_{q_{j_1}}) P_{q_{j_2}} \cdots P_{q_{j_k}} \Vert \le C 
 \frac{\vert q'_{j_1}-q_{j_1} \vert }{d(q_{j_1},X^{j_k -1})},
  \end{equation}
   for all  $1\le k \le \ell$. In particular, 
\begin{equation}\label{dequa3}
 \Vert P_{q'_{j_1}} - P_{q_{j_1}}\Vert \le C 
 \frac{\vert q'_{j_1}-q_{j_1} \vert }{d (q_{j_1},X^{j_1 -1})}.
  \end{equation}
\end{definition}

\begin{remark}
{\rm 
The inequality (\ref{dequa3}) says that if $q_{j_1}$ and $q'_{j_1}$ are sufficiently 
close points in  the stratum ${S}^{j_1}$, 
then the norm $\Vert P_{q'_{j_1}} - P_{q_{j_1}}\Vert$ is roughly 
given geometrically by the angle between the tangent spaces  $T_{q_{j_1}} {S}^{j_1}$ and 
$T_{q'_{j_1}} {S}^{j_1}$. More precisely, if 
$$\theta = \text{angle} (T_{q_{j_1}}{S}^{j_1},  T_{q'_{j_1}}{S}^{j_1})$$ 
then one has $\Vert P_{q'_{j_1}}- P_{q_{j_1}} \Vert \thickapprox \sin(\theta)$. 
This angle is controlled by the inequality (\ref{dequa3}).
Similarly one has $\Vert P^\perp_{q'_{j_1}}  P_{q_{j_1}} \Vert = \sin(\theta)$
(see also \cite[Figure 5]{BMT}).
\\
The inequalities (\ref{dequa1}) and (\ref{dequa2}) 
are more general in the sense that the angles of 
all the tangent spaces at points appearing in the chain are taken into account.
}
\end{remark}

In the following, one provides the relation between the different stratifications:
Lipschitz, Kuo-Verdier and Whitney.  
In fact one explains (or recalls the proof of) the relation
\begin{center} 
Lipschitz (L) $\Rightarrow$ Kuo-Verdier (w) $\Rightarrow$ Whitney (b) $\Rightarrow$ Whitney (a).
\end{center}

\begin{remark}
{\rm
a) The Lipschitz conditions are stronger than Whitney conditions \cite{Whitney1965}.
First, let us remark that for a chain of two points $\{q_{j}, q_{k} \}$, then the norm 
$\Vert P^\perp_{q_{j}} \, P_{q_{k}} \Vert$ in the inequality  (\ref{dequa1}) is 
geometrically given by $\sin(\theta)$ where $\theta$ is the angle between the two tangent spaces 
$T_{q_{j}} {S}^{j}$ and $T_{q_{k}} {S}^{k}$.

Now, recall that the Whitney condition $(a)$ for a pair of adjacent strata ${S}^{j}$ and 
${S}^{k}$ with ${S}^{k} \subset \overline{{S}^{j}}$, says that the angle between 
the tangent spaces $T_y {S}^{k} $ and $T_x {S}^{j} $ tends to zero when 
$x\in {S}^{j}$ tends to $y\in {S}^{k}$ (see \cite{Trotman2022,Whitney1965}).  

The Lipschitz conditions, showed  in the inequality (\ref{dequa1}), are stronger than Whitney condition (a) 
in the sense that this condition is satisfied for all pairs of strata involved in the considered chain. 

b) The Lipschitz conditions are stronger than Kuo-Verdier conditions (w) 
which are themselves stronger than  Whitney condition (b) \cite{Juniati2017,Verdier1976,Trotman2022}. 

The Kuo-Verdier conditions (w) says that, in a neighbourhood sufficiently small
of the stratum ${S}^{k}$ (in $\C^m$), the angle $\theta$ between the tangent spaces 
$T_x {S}^{j}$ and $T_y {S}^{k}$ is estimated by an inequality 
\begin{equation}\label{deKuo}
\sin(\theta) \le G \vert x-y \vert,
\end{equation} for some constant $G$.

The  Whitney condition (b) is verified  for a pair of adjacent strata ${S}^{j}$ and 
${S}^{k}$ with ${S}^{k} \subset \overline{{S}^{j}}$, if for each sequence 
$\{ x_i\}$ of points in ${S}^{j}$ and each sequence $\{ y_i\}$ of points in 
${S}^{k}$, both converging to the same point $y$ in ${S}^{k}$, 
such that the sequence of secant lines 
$(x_i,y_i)$ converges to a line $\ell$ and the sequence of tangent 
spaces $T_x {S}^{j}$ converges to a space $T$ (in a suitable Grassmannian space) 
as $i$ tends to infinity, then $\ell$ is contained in $T$.

Now apply the inequality  (\ref{dequa1}) for two adjacent strata $({S}^{j}, {S}^{k})$. 
We have 
$$\sin(\theta) \le C \frac{\vert x - y \vert}{d(x, X^{j-1})}.$$
In formula (\ref{deKuo}), put 
$$G = \frac{C}{ \inf_{x\in \U} \{d(x, X^{j-1}) \}}$$
where $\U$ is a sufficiently small neighbourhood of ${S}^k$ (in $\C^m$). Then the 
 Kuo-Verdier conditions (w) are satisfied.
 
 We emphasize again that the Lipschitz conditions are stronger than the Whitney condition (b) 
 in the sense that this condition is satisfied for all pairs of strata appearing in the chain.}
 
\end{remark}

\bigskip

\section{Lipschitz vector fields}\label{ideas}

In  \cite{Schwartz1997} Marie-H\'el\`ene Schwartz shows that 
considering stratified vector fields is not sufficient to prove the Poincar\'e-Hopf theorem
for singular stratified varieties.
She gives counterexamples
for which the index of the considered vector field at an isolated singular point calculated 
on the stratum is not the same when calculated in the ambient space (see also \cite[\S 2.1]{BMT}).
She shows that the theorem can be proved when considering radial extensions, 
obtained as the sum of two extensions: parallel and transversal. 
She constructs these extensions in the framework of Whitney stratifications, 
that is very technical and delicate. 
In \cite{Schwartz2000}, Marie-H\'el\`ene Schwartz 
 uses the same construction to define Chern classes of stratified varieties.

In \cite{BMT}, we show how the Lipschitz stratifications framework allows us 
to obtain a simpler proof of the Poincar\'e-Hopf theorem, in 
particular concerning  the construction of parallel and transversal extensions. 
This later is named ``basic normal'' in  
\cite[\S 4.2]{BMT}, we will keep the terminology. 
The aim of this section is to show that the construction of these extensions 
in the Lipschitz  framework allows us to simplify the definition of characteristic classes
of complex singular analytic varieties. 


The main advantage of using Lipschitz stratifications is that, instead of constructing the parallel and transverse extensions by pairs of strata, the construction is done all strata at once, taking into account the distance to the strata to determine the angle of the extension.
Moreover, the continuity of constructed vector fields (resp. $r$-frames)
     comes directly from Lipschitz properties without using any homotopy.

\subsection{The parallel Lipschitz extension} \label{deparallel}

Lipschitz stratifications have (among others) the following property 
(see Mostowski (1985) \cite[Prop. 2.1]{Mostowski1985}):

\begin{proposition}\label{propos}
Let $\{ X^j \}_{j=1}^n$  be a Lipschitz stratification of $X$ 
and let $v$ be a Lipschitz stratified vector field on $X^j$. 
Then $v$ can be extended to a Lipschitz stratified vector
 field on ${S}^{j+1}$.
\end{proposition}

A Lipschitz vector field $v$ defined on $X^j$ and tangent to all strata of dimension 
less than (or equal to) $j$ can be extended  to a Lipschitz vector field $\widehat v$ defined on $\C^m$ 
with the same Lipschitz constant (see Kirszbraun Theorem 1934, for instance in  \cite{Federer2014}). 
Mostowski defines the parallel extension of $v$ by
\begin{equation}\label{dede}
\tilde v(x) = \begin{cases}
\widehat v(x) & \text{ if } x \in X^j \cr
P_x \widehat v(x) & \text{ if } x\in  {S}^{j+1}. \cr
\end{cases}
\end{equation}
Then $\tilde v$ is a stratified Lipschitz vector field defined on ${S}^{j+1}$ 
(no matter which extension $\widehat v$ is chosen). 

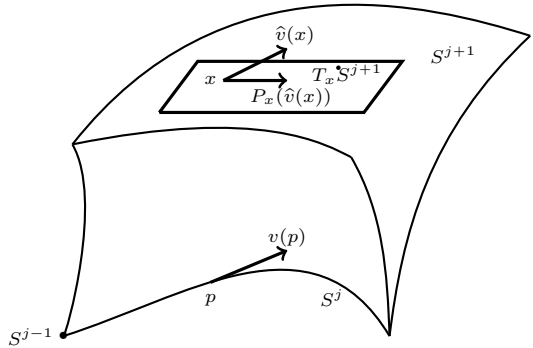
\begin{figure}[ht]\centering

\begin{tikzpicture}[scale=1.7]
  
\draw [thick](0.12,2) .. controls (1.3,3.5) and (2.9,3.1) .. (3.7,2.85);  
\draw [thick](0.12,2) .. controls (1.6,2.3) and (2.1,2) .. (2.3,1.9);  

\draw [thick](2.6,0.5) .. controls (2.7,1.3) and (2.9,2.1) .. (3.7,2.85);   
 \draw [thick] (2.6,0.5) .. controls +(92:0.5cm) and +(300:0.5cm) .. (2.3,1.9);   
 
 \draw [thick] (0.05,0.5) .. controls (0.07,0.5) and +(300:0.5cm) .. (0.12,2);
\draw [thick] (0.05,0.5) .. controls (1,0.8) and (2,1.5)  .. (2.6,0.5);
         
   \fill [black,opacity=1] (2.2,2.6) circle (0.5pt);
 \draw [very thick] (0.8,2.25) -- (1.1,2.65) -- (2.7,2.65) -- (2.4,2.25) --(0.8,2.25);
  
    \node at (2.27,2.55) {\scriptsize $T_x \strate^{\jalpha+1}$};
    
       \node at (1.2,2.5) {\scriptsize $x$};
\node at (1.86,2.85)  {\scriptsize ${\widehat v}(x)$}; 
\node at (1.83,2.37)  {\scriptsize $P_x({\widehat v}(x))$};     
 
\node at (1.8,1.27)  {\scriptsize $v(p)$};  
\draw[very thick,->] (1.3,2.5) -- (1.8,2.75);   
\draw[very thick,->] (1.3,2.5) -- (1.8,2.5);   

  \draw[very thick,->] (1.2,0.92) -- (1.8,1.17);  
\node at (1.2,0.9) [below] {\scriptsize $p$};  
    
\node at (2.15,0.8) {\scriptsize $\strate^{\jalpha}$}; 
\node at (3.1,2.7) {\scriptsize $ \strate^{\jalpha+1}$}; 
 
 \node at (0.05,0.5)  {\scriptsize $\bullet$}; 
  \node at (0.05,0.5) [left] {\scriptsize $ \strate^{\jalpha-1}$}; 

\end{tikzpicture}
\caption{The parallel Lipschitz extension.}\label{defig7}
\end{figure}

The following results are direct consequences of the construction
and of the definition of Lipschitz contidions
(angles between tangent spaces to the strata are small).

 \begin{corollary} \label{cor1}
 If $v_1$ and $v_2$ are $\C$-linearly independent on  a stratum ${S}^j$, 
then their Lipschitz parallel extensions are $\C$-linearly independent in a neighbourhood of  ${S}^j$.
\end{corollary}
 
 In a more general way, one has:
 \begin{corollary} \label{cor2}
If $v^{(r)}$ is an $r$-frame in a stratum ${S}^j$,
then its Lipschitz 
parallel extension is an $r$-frame in a neighbourhood of  ${S}^j$.
 \end{corollary}
 \medskip

 \subsection {The basic normal Lipschitz vector field}\label{S3.2}
\medskip

The M.-H. Schwartz' idea is given a field of vectors tangent to a stratum 
and having an isolated singular point, how to extend it in a neighbourhood 
of the singular point into a stratified vector field, with the same isolated singular point 
and the same index (calculated in the ambient space) at that point as the initial vector field.

The previous paragraph and the corollaries \ref{cor1} and \ref{cor2}  show how to extend 
an $r$-field tangent to a stratum 
into a stratified $r$-field in a neighbourhood of this stratum. 
This technique provides an  
$r$-field comprising a set of non-isolated singular points, namely a disk, 
slice of the tubular neighbourhood of the stratum.

In order to meet the conditions of the extension, the idea is to add 
to the preceding parallel $r$-field a ``transversal'' 
vector field ``tangent'' to the slice and having an isolated singular point
of index +1 in the center of the disk. This vector field is in a sense 
orthogonal  to the stratum, 
it is constructed in the real setting and named by ``basic normal'' vector field in 
\cite[\S 4.2]{BMT}.

The following Lemmas have been proved 
in the real case of semi-analytic sets in \cite{BMT}. 
The proof of the complex analytic case is similar.

Let $X\subset \C^m$ be a smooth analytic subset of some compact region, 
of dimension $n$ and let $X'  = \overline{X} \setminus X$. 
For $q\in \C^m$, let us denote by  $d_{X}(q) = d(q, X) $ 
 the distance from $q$ to $X$. 

\begin{lemma} \label{Lem1} \cite[Lemma 4.1]{BMT} 
There exists an integer 
$K \in \N$ and an $\varepsilon >0 $ such that for 
$$U = \{ x \in \C^m : {d}_X (x) < \varepsilon {d}_{X'}^K (x) \},$$
then for very $x\in U$ there exists exactly one point $\pi (x) \in X$ for which 
${d}_X(x) = \vert x - \pi(x) \vert$. The mapping $\pi$ is smooth in $U$. 
\end{lemma}

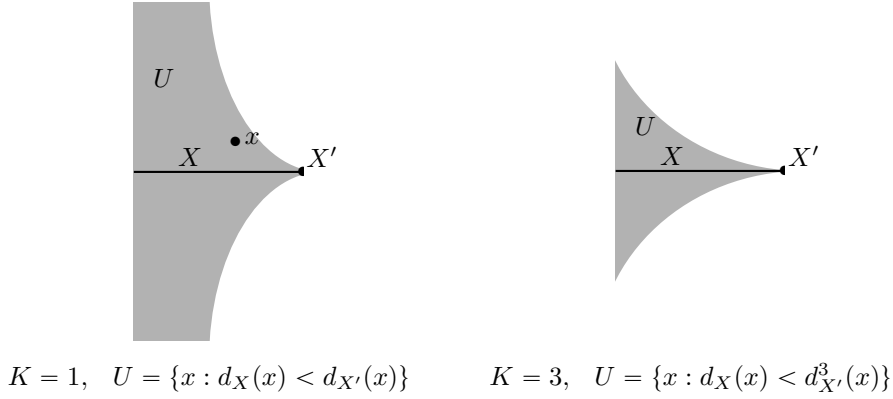
\begin{figure}[ht]\centering 
\begin{tikzpicture}[scale=0.5]

    \node at (7,0.4) {$X'$} ; 
  \node at (4,-5.5) {$K=1$, \; $U = \{ x : \arho_X(x) < \arho_{X'}(x) \}$} ;     
\clip (2,-4.5) rectangle (6.5,4.5);

 \fill [gray!60] rectangle (6.5,4.5);
 \fill [gray!60] rectangle (6.5,-4.5);
 
 \fill [white] (7,5) ellipse (3 and 5);
  \fill [white] (7,-5) ellipse (3 and 5);
  
   \node at (4.7,0.8) {$\bullet $} ; 
      \node at (4.7,0.9) [right] {$x$} ; 
   \draw[thick] (2,0) -- (6.5,0);
      \node at (2.8,2.5) {$U $} ; 
      \node at (3.5,0.4) {$X$} ; 
  \node at (6.5,0) {$\bullet $} ; 
  
\end{tikzpicture}
\qquad 
\begin{tikzpicture}[scale=0.5]

    \node at (7,0.4) {$X'$} ; 
  \node at (4,-5.5) {$K=3$, \; $U = \{ x :  \arho_X(x) < \arho^3_{X'}(x) \}$} ;     
\clip (2,-4.5) rectangle (6.5,4.5);

 \fill [gray!60] rectangle (6.5,4.5);
 \fill [gray!60] rectangle (6.5,-4.5);
 
 \fill [white] (7,5) ellipse (5.5 and 5);
  \fill [white] (7,-5) ellipse (5.5 and 5);
  
   \draw[thick] (2,0) -- (6.5,0);
   \node at (2.8,1.2) {$U $} ; 
      \node at (3.5,0.4) {$X$} ; 
  \node at (6.5,0) {$\bullet $} ; 
  
\end{tikzpicture}

\caption{The neighbourhood $U$.} \label{defig4}
\end{figure}

Define the {basic }  normal  vector field $r_X$ by $r_X(x) = x - \pi(x)$ after a translation to $x$. 

\begin{lemma} \label{Lem1b} \cite[Lemma 4.3]{BMT} 
The  {basic }  normal vector field  $r_X$ is smooth, Lipschitz in $U$ and 
satisfies $\Vert r_X(x)\Vert = {d}_X (x)$. 
It can be extended to a Lipschitz vector field on $\C^m$, again denoted by $r_X$, so that 
$\vert r_X \vert \le {d}_X$ and $r_X=0$ 
outside of $\widetilde U = \{ x : {d}_X(x)  < 2 \varepsilon {d}^K_{X'} (x) \}$.
\end{lemma}

In the following, we will apply these two Lemmas \ref{Lem1} and \ref{Lem1b} to 
 the situation of the stratum ${S}^j$ and will denote by $U_j$ and $r_{{S}^j} = r_j$ 
 respectively the neighbourhood and the {basic }  normal vector field corresponding.

 \subsection {``Champs sortants'' in Lipschitz stratifications}\label{desortant}
 
Let ${\mathcal S} = \{ X^j \}$ be a Lipschitz 
stratification of the analytic variety $X$, subset of some compact region in $\C^m$.

The notion of ``champ sortant'' (literally ``outgoing field'') has been introduced by Marie-H\'el\`ene 
Schwartz under the terminology ``transversal vector field '' in her proof of the Poincar\'e-Hopf theorem, using Whitney stratification. In the previous paper 
\cite[\S 4.2]{BMT} we introduce the
notion of ``champ sortant'' in the framework of Lipschitz stratifications of real varieties. 
This construction applies also in the complex setting. 
Here we give a summary of this construction with the important elements.

Let  $q\in {{S}}^{j}$. Consider the space $T^\perp_q {{S}}^{j}$ as a 
subspace of $U_j$ (it is $\pi^{-1} (q)$ in the previous terminology). 
If $x \in T^\perp_q {{S}}^{j}$ and $z \in T_x \C^m$,  let $\omega_j $ 
be the orthogonal projection of $z$ onto $T_x (T^\perp_q {{S}}^{j})$.

\begin{definition}\label{desortie} (see \cite[Definition 5.1]{BMT} in the real setting). 
Let $\lambda$ be a constant, $0< \lambda <1$. 
A $\lambda$-``champ sortant'' is a Lipschitz stratified vector field $w$ on $\C^m$ such that, 
for every $j$, in some neighbourhood $U'_j \subset U_j$ of ${{S}^j}$, 
$$\Vert \omega_j  (w- \psi_j r_j) \Vert \leq \lambda \Vert \psi_j r_j \Vert = 
\lambda \Vert \omega_j  (\psi_j r_j) \Vert $$
for some (scalar valued) functions $\psi_j $ taking values in some interval $[a_1, a_2]$ with 
$0< a_1< a_2 <  \infty$; and so
$$ \sangle (w, r_j) < \lambda \quad  \text{in } U'_j.$$
Actually, the function $\psi_j $ will be constant. 
\end{definition}

For each point $q\in X$ and each stratum $X^j$, let us denote by 
$d_j(q) = d_{X^j}(q)$ the distance from $q$ to $X^j$, and, for the points $q,q'$
$$d_j(q,q') = \min (d_j(q),d_j(q')).$$ 

Consider a pair of strata  $({S}^k, {S}^j)$
such that $S^j \subset \overline{S^k}=X^k$. 
By  \cite[Lemma 1.1]{BMT} which applies also in the complex setting, 
for $q \in {S}^j, \; q' \in {S}^k$, 
there exists numbers $C>0$, $\beta >0$ such that
$$ \sphericalangle (q'-q, T_{q}  {S}^j) \lesssim C \frac{ \vert q' -q \vert^\beta}{{d}_{j-1} (q,q')}.$$

It follows that if $\vert q'-q \vert \lesssim \min({d}_{j-1} (q))^K, \frac{1}{2} {d}_{j-1} (q))$, then
$$\sphericalangle (q'-q, T_{q}  {S}^j ) \le C \frac{ \vert q' -q \vert^\beta}{{d}_{j-1} (q)}
 \lesssim {d}_{j-1} (q)^{j K}.$$

It follows that if $K$ appearing in the definition of $U_j$ satisfies the estimate 
$K\ge 2 / \beta$, 
then in our stratification, we have the above estimate in  $U_j$. 

Let us fix a stratum ${{S}}^{j}$ and take $K$ so big that in
$U_j^* = \{ {d}_j < 2 {d}^K_{j-1} \}$
we have the  {basic }  normal vector field $r_{{S}^j}= r_j$ with all
properties of Lemma \ref{Lem1b}  
and the above form of (b) holds in $U_j^*$. Let 
$U_j = \{ {d}_j <{d}_{j-1}^K \}$, and $\varepsilon = j K$. 

The proof of the following Lemma (\cite[Lemma 5.2]{BMT} in the real setting)
applies in the complex setting:

\begin{lemma} \label{Lemma5.1}
There exists $\varepsilon>0 $ and 
a Lipschitz ${\mathcal S}$-stratified vector field $w_j$ such that
$$\Vert w_j - r_j \Vert \le {d}_j^{\varepsilon + 1} \text{ in } U_j,$$
$$\Vert w_j \Vert \le \Vert  r_j \Vert = {d}_j  \text{ in } U_j^*,$$
$$ w_j = 0 \text{ outside of } U_j^*.$$
\end{lemma} 


\begin{remark}{\rm 
Of course, we have $\Vert w_j\Vert \le \min ({d}_j, {d}_{j-1}^K)$. }
\end{remark}

\begin{proposition} (\cite[Proposition 5.4]{BMT} in the real setting). 
Let $w_i$'s be as in Lemma \ref{Lemma5.1} 
and $\lambda$ any number, $0<\lambda <1$. 
Then, if $0< c_0 < c_1< \cdots < c_n$ 
in a sufficiently rapidly increasing sequence of numbers 
({\it i.e.} $\frac{c_{i+1}}{c_i}$ is sufficiently big, for all $i$), 
then the vector field $w = \sum c_i w_i$ is a ${\mathcal S}$-stratified Lipschitz vector field 
 $\lambda$-``champ sortant''.
\end{proposition}

\begin{proof}
Let us fix a stratum ${{S}}^{j}$  and consider $w_j$ in the neighbourhood $U_j$ 
of ${{S}}^{j}$. 

Let $i<j$. Since $w_i$ is Lipschitz and stratified, we have, in $U_j$, 
$$\Vert \omega_j (w_i) \Vert \le {d}_j$$
because $\omega_j(w_i)=0$ on ${{S}}^{j}$. Since 
$$\Vert \omega_j (w_j) \Vert \sim \Vert w_j \Vert \sim  {d}_j,$$
we have, for a sufficiently big constant $A$,
$$\Vert \omega_j(w_i) \Vert < A \Vert \omega_j (w_j)\Vert.$$

Let $i>j$. Then $\Vert w_i \Vert \le {d}^K_{i-1} \le {d}^K_j$ (here $K$ is as in 
Lemma \ref{Lemma5.1}; we may assume that $K\ge 2)$. It follows that, no matter how big is 
a constant $A$, in a suitable neighbourhood of  ${{S}}^{j}$, we have 
$$\Vert \omega_j (w_j) \Vert > A \Vert \omega_j (w_i) \Vert.$$
A bound below for $\omega_j (w_j)$ is given by 
$$\Vert \omega_j (w_j)\Vert \sim \Vert w_j \Vert \sim d_j.  $$
Thus, if $c_i$ increases sufficiently rapidly,
$$ \Vert w - c_j w_j \Vert =\Vert  \sum_{i\ne j} c_i w_i \Vert < \frac{\lambda}{2} \, \Vert c_j w_j \Vert,$$
and, by Lemma \ref{Lemma5.1}
$$ \Vert w - c_j r_j \Vert < \lambda \, \Vert c_j r_j \Vert$$
in some neighbourhood of ${{S}}^{j}$.
\hfill
\end{proof}

The ``champ sortant'' is the stratified Lipschitz vector field $w$.

\medskip

 \subsection{Applications}

According to section \ref{dobstruction}, the aim of the construction is to attach to each $2p=2(m-r+1)$-cell 
of a cellular dual decomposition, the index of a complex $r$-field with isolated singularities. 
The $r$-field will be constructed by using the parallel Lipschitz extension
and the ``champ sortant''.

The parallel Lipschitz extension and the Lipschitz ``champ sortant'' are well defined
within a neighbourhood ${\mathcal U}_j$ for each stratum ${S}^j$. That means that we have 
to refine the triangulation $(K)$ and consequently also the dual cellular decomposition 
so that each dual cell meeting the stratum ${S}^j$  is included in 
${\mathcal U}_j$, this for all strata.

Now, let $d$ be a $2p$-cell dual of a $(K)$-simplex $\sigma^{2(r-1)} \subset {S}^j$,
where ${S}^j$ is a stratum of (complex) dimension $j \ge (r-1)$. 

The obstruction dimension to construct an $r$-frame tangent
in ${S}^j$ is $2p_j = 2(j-r+1)$, that is the dimension of  $d \cap {S}^j$
(note that $d$ is topologically transverse to $S^j$ - see Lemma \ref{delemme}).
One can construct a Lipschitz $r$-frame $v^{(r)} = (v_1,\ldots, v_{r-1}, v_r)$ on $d \cap {S}^j$
 with an isolated singularity due to cancellation  
 of the last vector $v_r$  at the barycenter $\hat d$ of the cell $d$. 
Let us recall that $\hat d$ is also barycenter of the simplex $\sigma^{2(r-1)}$.
 
 Let us denote by ${\overline v}^{(r)}$
 the parallel Lipschitz extension of $v^{(r)}$ (defined in section \ref{deparallel}) 
 to the last vector of which we add 
 the {``champ sortant''} $w$ (defined in section \ref{desortant}), that means
 $$ {\overline v}^{(r)} = (\widetilde v_1,\ldots, \widetilde v_{r-1}, \widetilde v_r + w).$$
 That is a Lipschitz stratified  $r$-field well defined on the $2p$-cell $d$ with an isolated 
 singularity at the barycenter $\hat d$ of  $d$.
 
 \begin{proposition} \label{gunha} 
 (a) Let $v^{(r)} = (v_1,\ldots, v_{r-1}, v_r)$ be a Lipschitz $r$-field tangent to ${S}^j$, 
 with $j >r-1 $ 
 and with an isolated singularity at $\hat d \in {S}^j$ due to cancellation  of the last vector $v_r$.  
 Let us denote by ${\overline v}^{(r)}$
 the parallel extension of $v^{(r)}$ to the last vector of which we add the 
Lipschitz champ sortant: 
 $${\overline v}^{(r)} = (\widetilde v_1,\ldots, \widetilde v_{r-1}, \widetilde v_r + w). $$
 The point $\hat d$ is also an isolated singularity of ${\overline v}^{(r)}$ and one has equality of indices
 $$I(v^{(r)}, \hat d ;{S}^j) = I({\overline v}^{(r)}, \hat d;\C^m)$$
where $I(v^{(r)}, \hat d ;{S}^j)$ is the index of the $r$-field $v^{(r)}$ 
at its singular point $\hat d$,  calculated on the stratum $S^j$,
and $I({\overline v}^{(r)}, \hat d;\C^m)$ 
is the index of its extension, the $r$-field ${\overline v}^{(r)}$, in the ambient space, 
at its same singular point $\hat d$,\\
\noindent (b)  If $j = r-1$, then one can consider only the $(r-1)$-field 
$(v_1,\ldots, v_{r-1})$ on ${S}^j$.
 Then the $r$-field $ (\widetilde v_1,\ldots, \widetilde v_{r-1},  w)$ is 
 a Lipschitz stratified  $r$-field well defined on $d$ and its index $I({\overline v}^{(r)}, \hat d;\C^m)$
 is equal to $+1$.
  \end{proposition}
 
\bigskip
   
\subsection{Radial extension of $r$-frames} 

In the same way as in (\cite{BS}, \cite{Schwartz1965}, \cite{Schwartz2000}), we prove the following theorem, 
existence of radial extension of frames, 
in the Lipschitz framework:

\begin{theorem} \label{teora}
Let  $X$ be an analytic variety, subset of some compact region 
 in the complex space $\C^m$  and $( {\mathcal S} )$ 
a Lipschitz stratification of $\C^m$ which is composed of a stratification of $X$ to which we 
 add the stratum $(\C^m \setminus X)$.
Let $(K)$ be a triangulation of  $\C^m$ compatible with the stratification and denote by 
$(D)$ a cellular decomposition of  $\C^m$ dual of $(K)$. 
Let $ \mathcal D$ be the cellular neighbourhood of $X$ in $\C^m$, 
union of cells which meet $X$. 
We can construct, on the $2p$-skeleton $(D)^{2p}$, 
a  Lipschitz stratified $r$-field $v^{(r)}$ which is called {\rm radial field}, whose singularities 
satisfy the following properties:\\
\indent  (i)  $v^{(r)}$ has only  isolated singular points, which are zeroes of the last
vector field $v_r$ \and which are barycenters $\hat d$ of $2p$-cells.
 On $(D)^{2p-1}$, the $r$-field $v^{(r)}$ has no singular point 
and on $(D)^{2p}$ the $(r-1)$-field $v^{(r-1)}$ has no singular point. \\
\indent  (ii) Let $\hat d \in {S}^j \cap (D)^{2p}$ be a singular point of 
$v^{(r)}$ in the $2j$-dimensional stratum ${S}^j$. If $j>r-1$, the index of
$v^{(r)}$ at $\hat d$, computed in the ambient space is the same as the index of the
restriction of $v^{(r)}$ to ${S}^j \cap (D)^{2p}$ considered as an $r$-field
tangent to ${S}^j$. If $j=r-1$, then 
 $I(v^{(r)}, \hat d) = +1$.\\
\indent  (iii) Inside a $2p$-cell $d$ which meets several strata, the only singularities
of $v^{(r)}$ are inside the lowest dimensional one (in fact located at the
barycenter of $d$). \\
\indent  (iv) The $r$-field $v^{(r)}$  has no singularity on $\partial \mathcal D \cap (D)^{2p}$.
\end{theorem}

Basically, the proof follows the same process as the proof of the Poincar\'e-Hopf theorem
\cite[\S 6]{BMT}. The proof of the Poincar\'e-Hopf Theorem 
 is a proof by induction starting with a radial vector field pointing outwards balls centered 
at points corresponding to the $0$-dimensional strata. Here the role of the strata of dimension $0$ 
is played by the strata ${S}^j$ such that 
the intersection of a cell $d$ of dimension $2p= 2(m-r+1)$ with ${S}^j$ 
is the barycenter  $\widehat d$ of $d$. We have therefore $2 (j -r+1) = 0$, i.e. $j= (r-1)$.

\begin{proof}
{Let ${S}^j$ be a stratum  of (real) dimension $2 j= 2(r-1)$.
The intersection of a cell $d$ of dimension $2p= 2(m-r+1)$ with ${S}^j$ 
is the barycenter  $\widehat d$ of $d$. }
At the points $\widehat d$, we give ourselves a complex Lipschitz $(r-1)$-field tangent to ${S}^j$.
It can be extended along $d$ by Lipschitz extension, with the condition that $d$ 
is sufficiently small to be included in $\U_j$. 

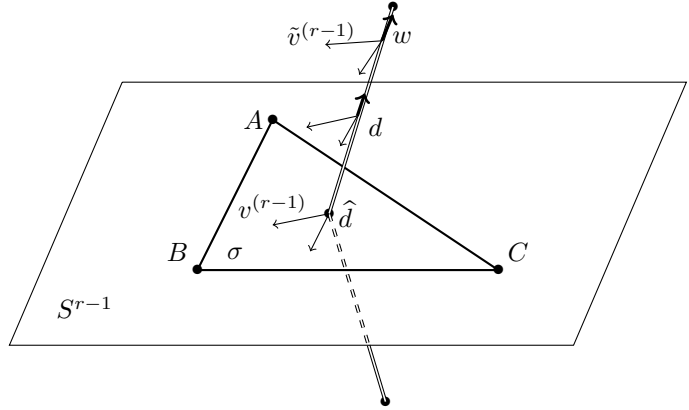
\begin{figure}[ht]\centering
\begin{tikzpicture}[scale=0.5]

\coordinate (A) at (-5,3);
\coordinate (B) at (-7,-1);
\coordinate (C) at (1,-1);
\draw [thick](B)--(A)-- (C) -- (B);
\node at (B)  {$\bullet$};   
\node at (A)  {$\bullet$};  
\node at (C)  {$\bullet$};    
\node at (B)[above left]  {$B$}; 
\node at (A)[left]  {$A$};   
\node at (C)[above right]  {$C$};   

\coordinate (X) at (-3.5,0.5);

\coordinate (G) at (-1.8,6);
\coordinate (K) at (-2,-4.5);
\node at (G)  {$\bullet$};   
\node at (K)  {$\bullet$};  

\node at (X)  {$\bullet$};  
\coordinate (S) at (-12,-3);
\coordinate (T) at (3,-3);

\coordinate (AE) at (intersection of X--K  and  S--T);

\draw[double](G) -- (X);
\draw [double,dashed](X) -- (AE);
\draw[double] (AE)--(K);

\node at (-6,-0.55)  {$\bf \sigma$};
\node at (-2.7,2.8) [right] {$d$};

\coordinate (R) at (-9,4);
\coordinate (S) at (-12,-3);
\coordinate (T) at (3,-3);
\coordinate (U) at (6,4);
\draw (R)--(S)--(T) -- (U) --(R);
\node at (-10,-2)  {$\strate^{r-1}$};   

\node at (X) [ right] {$\widehat d$};   

\draw [->]  (-3.5,0.5) --(-4,-0.5);
\draw [->]  (-3.5,0.5) --(-5,0.2);
\node at (-5,0.7) {${v}^{(r-1)}$};
\draw [->]  (-2.75,3.1) --(-3.2,2.3);
\draw [->]  (-2.75,3.1) --(-4.1,2.8);
\draw [->]  (-2.1,5.1) --(-2.7,4.2);
\draw [->]  (-2.1,5.1) --(-3.6,5);
\node at (-3.7,5.3){$\tilde{v}^{(r-1)}$};

\draw [very thick,->]  (-2.75,3.1)  --(-2.55,3.7);
\draw [very thick,->]  (-2.1,5.1) --(-1.8,5.8);
\node at (-1.55,5.2) {$w$};

 \end{tikzpicture}
  \caption{The starting point of induction (real picture of a complex situation). 
The cell $d$ is dual of the simplex $\sigma = ABC$.}\label{debut}
   \end{figure}
   
More precisely, the starting point of the induction is as follows: 

If ${S}^{r-1}$ is a stratum whose real dimension is $2r-2 = 2(m-p)$, the intersection of a 
$2p$-cell with ${S}^{r-1}$ is a point: the barycenter $\widehat d$ of $d$, which is also the 
barycenter of the simplex $\sigma$ 
(in ${S}^{r-1}$) of which $d$ is the dual cell (see Figure \ref{debut}).

The obstruction 
dimension for the construction of an $r$-frame tangent to ${S}^{r-1}$ is zero. 
One takes any Lipschitz  $(r-1)$-frame 
$v^{(r-1)}$ tangent to ${S}^{r-1}$ at the barycenter $\widehat d$ and the last vector 
$v_r$ is taken as zero at these points.

We construct the Lipschitz parallel extension 
(see section \ref{deparallel}) of the $(r-1)$-frame $v^{(r-1)}$ in the 
 cellular tubular neighbourhood ${\mathcal U}_{r-1}$
as an $(r-1)$-frame denoted by $\tilde{v}^{(r-1)}$.  
The extension, together with the Lipschitz  champ sortant $w$ 
 (see section \ref{S3.2}) provide an $r$-field
$\tilde{v}^{(r)} = (\tilde{v}^{(r-1)},w)$ which has  an  isolated singularity at $\widehat d$. 
According to Proposition \ref{gunha} (b), one has 
$I(\tilde{v}^{(r)},\widehat d)=+1$. 
\bigskip

Now, the induction process goes as follows: 
Let us suppose $k>r-1$ and the construction already performed on all strata ${S}^j$ 
whose (real) dimension is less than $2k$. That means that the construction provides 
an $r$-field on the $2p$-skeleton of a cellular neighbourhood $\D^j$ 
of ${S}^{j}$ for all strata ${S}^j$ with dimension $j < k$.
As in Proposition \ref{gunha}, 
 the $r$-field on ${S}^{j}$ is denoted by $v^{(r)} = (v^{(r-1)},v_r)$ 
and its Lipschitz extension by ${\overline v}^{(r)}$. 
If a $2p$ cell is transverse to a stratum ${S}^j$, then the $r$-field has an isolated 
singularity at the barycenter $\widehat d$, due to cancellation of the last vector $v_r$. 
The indices satisfy  then 
\begin{equation} \label{derechef}
I(v^{(r)}, \hat d ;{S}^j) = I({\overline v}^{(r)}, \hat d;\C^m).
\end{equation}
If the $2p$-cell $d$ does not meet a stratum and lies in the boundary of the 
cellular neighbourhood $\D^j$ of a stratum ${S}^j$, then the 
basic normal vector field $w$ does not vanish on $d$ and the $r$-field $\tilde{v}^{(r)}$ 
has no singularity on $d$. One remarks that such cells are dual of simplices in $(K)$ 
which have a face in ${S}^j$ but are not included in ${S}^j$. 

The boundary of ${S}^k$ is union of strata ${S}^j$, each of them 
equipped with a cellular neighbourhood $\D^j$ and an $r$-field 
satisfying (\ref{derechef}) at each singular point. 

The $r$-frame is constructed on the $2p$-skeleton of 
a cellular neighbourhood of the boundary of ${S}^k$. 

We extend the $(r-1)$-frame $\tilde v^{(r-1)}$ inside ${S}^k$, more precisely on the 
$2p_k$-skeleton of the intersection $(D)^{(2p)} \cap {S}^k$, 
with $2p_k=2(k-r+1)$. This is performed without singularity. 

Using the Lipschitz extension (\ref{deparallel}), one extends the last vector $v_r$ with 
singularities at the barycenters $\widehat d$ of cells $d$ in $(D)^{(2p)} \cap {S}^k$.
One obtains, on the stratum ${S}^k$ an $r$-frame with isolated singularities at 
the barycenters $\widehat d$, 
which are zeroes of the last vector $v_r$. One has (Proposition  \ref{gunha} (a))
$$I(\widetilde v^{(r)},\hat d; \C^m) = I(v^{(r)}, \hat d;{S}^k).$$

In summary, an $r$-frame already known on a neighbourhood of the boundary of a
stratum is extended with isolated singularities inside (a suitable
skeleton of) the stratum and then extended with property (a) of the Proposition \ref{gunha} to a
cellular neighbourhood of this stratum.

We denote by  ${\mathcal D}$ the tubular neighbourhood of $X$ in $\C^m$ consisting of
the $(D)$-cells which meet $X$. By subdividing if necessary, we can assume that ${\mathcal D}$ 
is contained in the union of the Lipschitz neighbourhood ${\U}^j$. 
\end{proof}

\section{Obstruction cocycles and classes}\label{declasse}

We define a $2p$-dimensional $(D)$-cochain in
$C^{2p}({\mathcal D},\partial {\mathcal D})$ by: 
\begin{equation}\label{classeS}
\gamma^p : d\quad  \to  \quad  I({\widetilde v}^{(r)},\widehat d).
 \end{equation}
The cochain $\gamma^p$ actually is a cocycle whose class $c^p(X)$ lies in   
$$H^{2p}({\mathcal D},\partial {\mathcal D})
\cong H^{2p}({\mathcal D},{\mathcal D}\setminus X) \cong H^{2p}(\C^m ,\C^m \setminus X),$$
where the first isomorphism is given by retraction along the rays of ${\mathcal D}$ 
 and the second by excision (by $\C^m \setminus {\mathcal D}$). 
 The class $c^p(X)$ coincides with the Schwartz class, because on the one hand, 
 the Lipschitz frame is a particular 
 case of the frames defined by M.-H. Schwartz and we know that the Schwartz class
 does not depend on the choice of frame in her framework.
On the other hand, a Lipschitz stratification is particularly a Whitney stratification.
 
\medskip
\begin{definition}
The $p$-th Schwartz class of $X$ is the  class 
 $$c^p(X)\in H^{2p}(\C^m ,\C^m \setminus X).$$
 \end{definition} 

Marie-H\'el\`ene Schwartz proved ``by hand'' that the resulting class is independent of the 
stratification, of the triangulation (compatible  with the stratification) and of the choice of the
radial $r$-frame. 
Here, within the Lipschitz framework, the resulting class is independent of the 
Lipschitz stratification (including the Lipschitz constant $C$), of the triangulation
 and of the choice of the Lipschitz extensions of frames. 
The proof of independence is now facilitated as we 
know \cite{BS} that by Alexander duality  isomorphism 
$$H^{2(m-r+1)}( \C^m , \C^m \setminus X) \to H_{2(r-1)} (X), \qquad d \to \sigma$$
image of the Schwartz class is the MacPherson class  \cite{MacPherson1974} 
of the constructible function  ${\bf 1}_X$ and that the MacPherson class 
does not use neither triangulation, nor $r$-frames. It is proved independent of the used stratification.  Nowadays, the classes are called Schwartz-MacPherson classes.


\begin{thebibliography}{99}

\bibitem {Brasselet1} J.-P. Brasselet, {\it An introduction to Characteristic Classes}, Publications IMPA, $33^o$  Col\'oquio Brasileiro de Matem\'atica. 

\bibitem {Brasselet-Thuy} J.-P. Brasselet and {\fontencoding{T5}\selectfont Th\h{u}y 
Nguy\~\ecircumflex n} T.B, 
{\it O Teorema de Poincar\'e-Hopf}, C.Q.D.- Revista Electr\^onica Paulista de Matem\'atica 
16 (2019), p. 134-162 (in Portuguese). 

\bibitem {BMT} J.-P. Brasselet, T. Mostowski and {\fontencoding{T5}\selectfont Th\h{u}y Nguy\~\ecircumflex n} T.B,  {\it Poincar\'e-Hopf Theorem for singular  analytic varieties. The Marie-H\'el\`ene Schwartz' ideas in the Lipschitz framework.} 
Banach Center Publications, Volume 128, Warszawa, 2024, 23--43.

\bibitem {BS} J.P. Brasselet and M.-H. Schwartz: {\em Sur les classes de
Chern d'un ensemble analytique complexe}, Ast\'erisque 82-83 (1981), 93-147.

\bibitem {Ch2} S.S. Chern, {\it Characteristic classes of hermitian 
manifolds}, Ann. Math. {\bf 47} (1946), 85-121.

\bibitem {DK} J. F. Davis and P. Kirk, {\it Lecture notes in algebraic topology}, Vol. 35. Graduate Studies in Mathematics. American Mathematical Society, Providence, RI, (2001).

\bibitem{Federer2014} H. Federer,  {\it Geometric measure theory}.  Classics in Mathematics, Springer 1996. 

\bibitem{Grothendieck2022} A. Grothendieck, {\it R\'ecoltes et Semailles, 
R\'eflexions et t\'emoignage sur un pass\'e 
de math\'ematicien, 1986},  Gallimard, 2023, available in english translation at 
\url{https://web.ma.utexas.edu/users/slaoui/notes/recoltes_et_semailles.pdf}

\bibitem{Juniati2017} 
D. Juniati, L. Noirel and D. Trotman,  {\it Whitney,  Kuo-Verdier and Lipschitz stratifications for the surfaces  $y^a = z^b x^c + x^d$.} Topology and  Applications, 234 (2018), 335-347.

\bibitem{Lefschetz1930} S. Lefschetz, {\it Topology}, Amer. Math. Soc, Coll. Publications, New York, 1930.

\bibitem{Lojasiewicz1964} S. \L ojasiewicz, {\it Triangulation of semi-analytic sets}, 
Ann. Scuola Norm. Sup. di Pisa,  Ser. 3, 18-4 (1964), 449-474.

\bibitem{Lojasiewicz1995} S. \L ojasiewicz, {\it On semi-analytic and subanalytic geometry}, Panoramas of 
Mathematics, Banach Center Publications, Volume 34, 
Institute of Mathematics, Polish Academy of Sciences, Warszawa 1995, pp 89 -- 104.

\bibitem{MacPherson1974}  R. D. MacPherson,  {\it Chern classes for singular algebraic varieties.} 
Ann. of Math. (2) 100, (1974), 423--432.

\bibitem{Mostowski1985}  T. Mostowski,  {\it Lipschitz Equisingularity}, Dissertationes Math. 
(Rozprawy Mat.) 243, PWN, Warszawa, 1985.

\bibitem{Mostowski1988} T. Mostowski,  {\it Tangent Cones and Lipschitz Stratifications}, Singularities, 
Banach Center Publications (S. Lojasiewicz, ed.), Vol. XX, PWN, Warszawa, 1988, pp. 303-322.

\bibitem{Munkres84} J. R. Munkres,  {\it Elements of algebraic topology}.
 Addison-Wesley Publishing Company, Menlo Park, CA, 1984.

\bibitem{Parusinski1988}  A. Parusi\'nski,  {\it Lipschitz Properties of Semi-analytic Sets}, Ann. Inst. Fourier, Grenoble, Vol. 38, 1988, no. 4, pp. 189-213.

\bibitem{Parusinski1993}  A. Parusi\'nski,  {\it Lipschitz stratifications}, in: Global Analysis in 
Modern Mathematics (Orono, 1991; Waltham, 1992), K. Uhlenbeck (ed.), Publish or Perish, 
Houston, 1993, 73--89.

\bibitem{Parusinski1994}  A. Parusi\'nski,  {\it Lipschitz stratification of subanalytic sets}, 
Ann. Sci. \'Ecole Norm. Sup. (4) 27 (1994), no. 6, 661--696. 

\bibitem  {Schwartz1965}  M.H.Schwartz, {\it Classes caract\'eristiques  d\'efinies par une stratification d'une 
vari\'et\'e analytique complexe,} C. R. Acad. Sci. Paris S\'er. I Math.   {\bf 260}, (1965), 3262-3264 et 3535--3537.

\bibitem  {Schwartz1997} M.-H.~Schwartz: {\em Champs radiaux sur  une stratification
analytique}, Travaux en cours, 39 (1997), Hermann, Paris.

\bibitem {Schwartz2000} M.-H. Schwartz, {\it Classes de Chern des ensembles
analytiques.} Actualit\'es math\'ematiques, Hermann, Paris, 2000.    

\bibitem {Siebenmann1979} L. Siebenmann and D. Sullivan,  {\it On Complexes that are 
Lipschitz Manifolds}, 
Geometric Topology (J. Cantrell, ed.), Academic Press, New York, 1979, pp. 503-552.

\bibitem {Steenrod1951} N. Steenrod,  {\it The Topology of Fibre Bundles} (1951). Vol. 14. Princeton Mathematical Series. Princeton University Press, Princeton, N. J.

 \bibitem {Trotman2022} D. Trotman, {\it Stratification theory,} Chapter 4  In: 
 Cisneros-Molina, J.L., L\^e, D.T., Seade, J. (eds.) 
Handbook of  Geometry and Topology of Singularities, Volume I. Springer, Cham, 231 -- 260.

 \bibitem {Verdier1976} 
 J. L. Verdier,  {\it Stratifications de Whitney et th\'eor\`eme de Bertini-Sard},  
 Invent. Math., Vol. 36, 1976, pp. 295-312.
 
 \bibitem {Whitney1965} 
  H. Whitney,   {\it Local properties of analytic varieties}, in Differential and Combinatorial Topology (ed. S. S. Cairns), Princeton Univ. Press (1965), 205-244.
 

\end{thebibliography}
\end{document}